\documentclass{amsart}

\usepackage[active]{srcltx}
\usepackage{euscript}
  {}%         % issue number
  {}%         % month
\newtheorem{theorem}{Theorem}
\newtheorem{lemma}[theorem]{Lemma}

\newtheorem{prop}[theorem]{Proposition}

\newtheorem{maintheorema}{Main Theorem A}
\newtheorem{maintheoremb}{Main Theorem B}

\theoremstyle{definition}

\theoremstyle{remark}
\newtheorem{remark}[theorem]{Remark}

\newcommand{\abs}[1]{\left\vert#1\right\vert}
\newcommand{\norm}[1]{\left\Vert#1\right\Vert}
\newcommand{\R}{\mathbb R}

\newcommand{\Lie}{\EuScript L}

\newcommand{\nablash}{\nabla{\kern -.75 em
     \raise 1.5 true pt\hbox{{\bf/}}}\kern +.1 em}
\newcommand{\Deltash}{\Delta{\kern -.69 em
     \raise .2 true pt\hbox{{\bf/}}}\kern +.1 em}
\newcommand{\Rslash}{R{\kern -.60 em
     \raise 1.5 true pt\hbox{{\bf/}}}\kern +.1 em}

\newcommand{\Vol}{\operatorname{Vol}}

\newcommand{\chib}{\bar\chi}

\newcommand{\Hb}{\bar H}

\newcommand{\Sphere}{\mathbb S}
\newcommand{\D}{\partial}

\title{Blow-up in the parabolic scalar curvature equation}
\author{Brian Smith}

\begin{document}

\address{
  Freie Universit\"at Berlin, Arnimallee 3, 
  14195 Berlin, Germany}

\email{bsmith@math.fu-berlin.de}

 %53C21 Methods of Riemannian Geometry 
 %35C44 Geometric Evolution Equations
 %35K55 Nonlinear PDE of Parabolic Type
 %35K57 Reaction-diffusion equations 
%\date{April 2006-\today.}


\begin{abstract}
The \textit{parabolic scalar curvature equation} is a
reaction-diffusion type equation on an $(n-1)$-manifold $\Sigma$,
the time variable of which shall be denoted by $r$.  Given a
function $R$ on $[r_0,r_1)\times\Sigma$ and a family of metrics
$\gamma(r)$ on $\Sigma$, when the coefficients of this equation
are appropriately defined in terms of $\gamma$ and $R$, positive
solutions give metrics of prescribed scalar curvature $R$ on
$[r_0,r_1)\times\Sigma$ in the form
\[
   g=u^2dr^2+r^2\gamma.
\]
If the area element of $r^2\gamma$ is expanding for increasing
$r$, then the equation is  parabolic, and the
basic existence problem is to take positive initial data at some
$r=r_0$ and solve for $u$ on the maximal interval of existence,
which above was implicitly assumed to be $I=[r_0,r_1)$; one often
hopes that $r_1=\infty$.  However, the case of greatest physical
interest, $R>0$, often leads to blow-up in finite time so that
$r_1<\infty$.  It is the purpose of the present work to
investigate the situation in which the blow-up nonetheless occurs
in such a way that $g$ is continuously extendible to $\bar
M=[r_0,r_1]\times\Sigma$ as a manifold with totally geodesic outer
boundary at $r=r_1$.
\end{abstract}

\maketitle

\section{introduction}

Given a smooth family of
Riemannian metrics $\gamma(r),r\in [r_0,\infty)$ on an
$(n-1)$-manifold $\Sigma$, the \textit{parabolic scalar curvature
equation} refers to the equation
\begin{equation}\label{eq:primus}
    \Hb r\frac{\D u}{\D r}=u^2\Delta_{\gamma}u+Au-\frac{1}{2}\left(\bar R-r^2R\right)u^3,
\end{equation}
where $\bar  R_r$ is the scalar curvature of $\gamma(r)$, the
function $R$ is ``arbitrary'', and the remaining terms in the
coefficients are defined by
\begin{gather*}
   A=r\frac{\D\Hb}{\D r} - \Hb +\frac12 \abs{\chib}_\gamma^2 + \frac12 \Hb^2, \\
   \chib_{AB}=\gamma_{AB}+\frac 12 r\frac{\D\gamma_{AB}}{\D r}, \\
   \Hb=\text{tr}_{\gamma}\chib=(n-1)+\frac 12 r\frac{\D\gamma_{AB}}{\D r}\gamma^{AB};
\end{gather*}
$A,B$ are used to denote components with respect to local
coordinates $\theta^i$ on $\Sigma$. Positive
solutions of Equation~\eqref{eq:primus} on an interval $[r_0,r_1)$
give metrics $g$ of prescribed scalar curvature $R$ on
$M=[r_0,r_1)\times\Sigma$ in the form
\[
    g=u^2dr^2+r^2\gamma.
\]
For more on  Equation~\eqref{eq:primus}, including  derivations, see~\cite{bartnik93},~\cite{SW},~\cite{SW2},~\cite{ST},~\cite{emw}.
The only derivation (that the author is aware of) in the $n$-dimensional case in the present context appears in~\cite{emw}.
This is done very quickly in that work, and so another derivation is provided in Appendix B.    
The function $\Hb$ and tensor $\chib$ are closely related to the extrinsic 
geometry of the hypersurfaces $\Sigma_r=\{r\}\times \Sigma$. Indeed, with $H,\chi$ the mean curvature and second fundamental form
of  $\Sigma_r$, one
has
\begin{align*}
      \Hb&=ruH\\
      \chib&=\frac{u}{r}\chi.
\end{align*}

In the case that $f\equiv r^2R/2-\bar R/2$ is positive and bounded
away from 0, it is easily established by using the maximum
principle that solutions will not exist for all $r>0$, but will in
fact blow up for some finite value of $r$.  It is the purpose of
the present work to investigate the blow-up behavior in the case that $\Sigma$ is compact.

The simplest case of blow-up, which we shall refer to as the
trivial case, occurs under the assumption that $f$ is fixed and
positive and $u$ is constant on each $\Sigma_r$. Then $\chib =
\gamma$ and $\Hb = (n-1)$, so that Equation~\eqref{eq:primus} is
reduced to the ordinary differential equation
\[
    (n-1)r\frac{d u}{dr}=\frac{(n-1)(n-2)}{2}u+fu^3.
\]
  For ``initial" data
$u(r_0)=u_0$, the solution of this problem is
\[
      u(r)=\frac{1}{\sqrt{c_0\left(\left(\frac{r_1}{r}\right)^{n-2}-1\right)}},
\]
where
$r_1^{n-2}=\frac{(n-1)(n-2)}{2f_0}u_0^{-2}r_0^{n-2}+r_0^{n-2}$
and $c_0=\frac{2f_0}{(n-1)(n-2)}$. Although the solution clearly
blows up at $r=r_1$ the metric $g=u^2dr^2+r^2\gamma$ is defined up
to and including $r_1$ as a $C^{\infty}$ metric on a manifold with totally geodesic outer boundary
boundary at $r=r_1$.  This is seen by making the
change of variables $\tilde r=\tilde r_0+\int_{r_0}^{r_1}udr$,
which puts the metric in the form $g=d\tilde r^2+r^2\gamma$.  It
is natural to ask: \textit{more generally, when can we expect this
behavior?} As a partial answer to this question, in this work the
following theorem is proved:
\begin{maintheorema}
Let $\Sigma$ be a compact $(n-1)$-manifold with a fixed metric $\gamma$. Let
$R$ be a $C^{\infty}$ function on
$[r_0,r_1+\varepsilon),\,\varepsilon>0$  such that $r^2R$ is
 non-decreasing and $f\equiv r^2R/2-\bar R/2> 0$.  Let $u$ be a  solution of
Equation~\eqref{eq:primus} on $[r_0,r_1)$ such that
\[
    \inf_{[r_0,r_1)\times\Sigma}u\sqrt{\left(r_1/r\right)^{n-2}-1}\geq\mu>0.
\]
Then $\lim_{r\to r_1}u\sqrt{\left(r_1/r\right)^{n-2}-1}$ exists and is a positive
$C^{\infty}$ function $\omega$ on $\Sigma$ so that the metric
\[
    g=u^2dr^2+r^2\gamma=\frac{\omega^2}{\left(r_1/r\right)^{n-2}-1}dr^2+r^2\gamma
\]
is extendable to $\bar M=[r_0,r_1]\times\Sigma $ in the sense
that $g\in C^{\infty}(M)\cap C^0(\bar M)$.
\end{maintheorema}
Thus, when $f$ is positive and non-decreasing, one can assert that
if $u $ blows up everywhere on $\Sigma$ at the blow-up time
$r=r_1$, and the blow-up happens \textit{at least} as fast as in
the trivial case, then the solution blows up exactly at this rate,
which after a change of variables allows the corresponding metric
to be extended in the sense of $C^0$ to the boundary component
$r=r_1$. Although the metric is not verified in this work to be $C^1$ at $r=r_1$, 
one can nonetheless compute the second fundamental form of $\Sigma_{r_1}$,
which is found to vanish.  That is, the outer boundary is totally geodesic.  

In order to ensure that the theorem is not trivially satisfied only, it 
is important to have examples of nontrivial blow-up.  
In Appendix~A non-trivial blow-up solutions in the case that 
$(\Sigma,\gamma)$ is the flat 2-torus are obtained from solutions 
of the curve shortening flow.  For non-trivial blow-up in the case $\Sigma=\Sphere^2$, see~\cite{FHS}.  
In that work the authors obtain
non-trivial blow-up using bifurcation theory and dynamical systems techniques.

The proof of the main theorem proceeds after the observation that
when $\gamma$ is fixed in $r$ some simple changes of variables
transform Equation~\eqref{eq:primus} into a more manageable form.
To see this,  note that when $\gamma$ is fixed $\chib = \gamma,\Hb
= n-1$ so that $A=(n-1)(n-2)/2$ and Equation~\eqref{eq:primus} is
\[
     (n-1) r\D_ru=u^2\Delta u+\frac{(n-1)(n-2)}{2}u+f u^3,
\]
where the subscript $\gamma$ on the Laplacian has been dropped as
will be done in the remainder.  Then  the function $\tilde u\equiv
r^{1-\frac{n}{2}}u$ verifies
\[
   (n-1)r^{(3-n)}\frac{\D\tilde u}{\D r}=\tilde u^2\Delta\tilde u+f\tilde
   u^3.
\]
Thus, defining
\[
    t=\frac{r^{n-2}}{(n-1)(n-2)},
\]
and regarding $\tilde u=\tilde u(p,t)$,
our equation takes the much nicer form
\begin{equation}\label{eq:secundus}
    \frac{\D\tilde u}{\D t}=\tilde u^2\Delta\tilde u+f\tilde
   u^3,
\end{equation}
Note that this equation has the scaling property that if $\tilde
u(p,t)$ is a solution then $\lambda \tilde u(p,\lambda^2 t)$ is
also a solution; this will be used below to assume without loss of
generality that the blow-up time occurs at $t=1$. Note also that
the order of the `space' and `time' variables has been switched
from what it was previously to the more standard order for
parabolic equations.

In the case that $f$  is  also fixed, examples of blow-up that occur
exactly like the special case discussed above can now be be
generated, in principle, by separation of variables: the function
$\tilde u=v/\sqrt{t_1-t}$ verifies Equation~\eqref{eq:secundus}
provided $v\in C^{\infty}(\Sigma)$ is a positive solution of
the \textit{stationary equation}
\begin{equation}\label{eq:stat}
      \Delta v+f_{t_1}v-\frac{1}{2v}=0.
\end{equation}
Following terminology as for the porous medium equation, solutions
$v(p)/\sqrt{t_1-t}$ generated in this way will be called
\textit{self-similar}. These solutions, if they exist, are very
special. But,  Main Theorem A asserts that in general if a
solution blows up at least as fast as the rate suggested by
the self similar blow-up, then in fact it blows up 
like a self-similar solution.

To prove that more generally blow-up is essentially self similar,
one can follow the same procedure used to generate self similar
solutions, with the generalization that the scaled function $v$ is
now allowed to depend on $t$. That is, defining
$v=\sqrt{t_1-t}\tilde u$, study the equation for $v$:
\[
     (t_1-t)\frac{\D v}{\D t}+\frac{1}{2}v=v^2\Delta v+fv^3
\]
Assuming without loss of generality that the blow-up  occurs at
$t_1=1$, a final change of  variables $t=1-e^{-\tau}$
yields
\begin{equation}\label{eq:veq}
      \frac{\D v}{\D\tau}=v^2\Delta v+fv^3-\frac{1}{2}v,
\end{equation}
and the blow-up behavior of the original equation can be dealt
with by studying the  behavior of $v$ as $\tau\to\infty$.
Specifically, the main theorem now follows from:
\begin{maintheoremb}
Suppose that $f>0$ is a $C^{\infty}$ function on
$[\tau_0,\infty)\times\Sigma$ such that for all $k\in\mathbb N$
\begin{align}
      \frac{\D f}{\D\tau}&\geq 0\label{cond:1}\\
       \norm{\left(e^{\tau}\D_{\tau}\right)^i(f-f_{t_1})}_{C^k(\Sigma)}&\leq C_k,\,\, i=0,1,2\label{cond:2'}
\end{align}
for constants $C_k$.
Let $v$ be a solution of Equation~\eqref{eq:veq}
on $[\tau_0,\infty)$ that satisfies
\begin{equation}\label{cond:2}
         v\geq\mu
\end{equation}
for some positive constant $\mu$. Then there exists a positive
solution $\omega\in C^{\infty}(\Sigma)$ of the stationary
equation, Equation~\eqref{eq:stat}, such that
$\lim_{\tau\to\infty}v=\omega$ in the sense of $C^k$ for any
$k\in\mathbb N$.
\end{maintheoremb}

The proof of this theorem, in turn, results from the successive
application of the next three theorems that will be proved in
the remainder.  
\begin{theorem}\label{thm:bound}
Assume $f,\,\D_{\tau}f\geq 0$.
Any  solution $v$ of Equation~\eqref{eq:veq} on an interval
$[\tau_0,\infty)$ satisfying $v\geq\mu$ for some positive constant
$\mu$ in addition satisfies $v\leq M$ for some constant
$M\leq\infty$.
\end{theorem}

\begin{theorem}\label{thm:profile}
Assume $\D_{\tau}f\geq 0$.
Let $v$ be  a  solution of Equation~\eqref{eq:veq}
on an interval $[\tau_0,\infty)$, which satisfies $\mu\leq v\leq M$ for
some positive constants $\mu,M$.  Then there exists a sequence
$\tau_i$ such that $v(\tau_i)$ converges uniformly to a positive
$C^{\infty}$ solution $\omega$ of the stationary equation.
\end{theorem}

Solutions of Equation~\eqref{eq:stat} will be referred to as
\textit{stationary states}.  This theorem asserts that the
$\omega$-limit set of $v$ is non-empty; it contains a stationary
state. As a consequence, we see that if the hypotheses of the
theorem are satisfied and in addition $f$ is fixed, then there is a
self-similar solution $\omega/\sqrt{1-t}$.

\begin{theorem}\label{thm:unique}
Assume Conditions~\eqref{cond:1}~and~\eqref{cond:2'} on $f$.
Let $v$ be a solution of~\eqref{eq:veq} on $[\tau_0,\infty)$
satisfying $\mu\leq v\leq M$, and let $\omega$ be a positive
$C^{\infty}$ stationary state in the $\omega$-limit set of $v$,
where convergence is taken in  the sense of $C^0$. Then $\omega$
is unique and $\lim_{\tau\to\infty}v(\tau)=\omega$, where the
limit can be taken in the sense of $C^k(\Sigma)$ for any $k$.
\end{theorem}
%Using standard parabolic theory, the convergence in the hypothesis
%can be taken sense of $L^p, p\geq 1$, but the previous theorem as
%stated suffices for present purposes.

The outline of the paper is as follows:

Section~\ref{sec:bounds} presents some basic pointwise
inequalities that are fundamental for most of the bounds  in the
remainder of the paper.  These inequalities are similar to
inequalities derived for the porous medium equation, originally by
Aronson-B\'{e}nilan~\cite{AB}.  The condition that $\gamma$ be fixed
is crucial.

Section~\ref{sec:global} is devoted to the proof of
Theorem~\ref{thm:bound}.   This is accomplished by proving a
strong global Harnack inequality for $v$ that shows that
$\sup_{\Sigma} v(\tau,p)$ is bounded in terms of $\inf_{\Sigma}
v(\tau+h,p)$ followed by a maximum principle argument that shows
that $\inf_{\Sigma} v(p,\tau+h)$ is globally bounded from above.

Theorem~~\ref{thm:profile} is proved in Section~\ref{sec:profile}.
This is done using  techniques similar to those used by
C.~Cortazar, M.~Pino, and M.~Elgueta in
~\cite{cpe1},~\cite{cpe2},~\cite{cpe3} to study blow-up in the
porous medium equation with source. The main tool is the
functional
\[
     J(v)=\int_{\Sigma}\left(|\nabla v|_{\gamma}^2-fv^2+\log v\right)dV_{\gamma},
\]
which is non-increasing  by virtue of Equation~\eqref{eq:veq} and
Condition~\eqref{cond:1}.  The bounds $\mu\leq v\leq M$ then show
that $J(v)$ is bounded from below, which leads to the existence of
a sequence $\tau_n$ such that $v(\cdot,\tau_n)$ converges to a
stationary state weakly in $H^1$ and strongly in $L^1$.

%Uniform convergence is obtained from $L^1$ convergence via results
%from Section~\ref{sec:local}.

%A marked difference between this procedure in the present case and
%the case of the porous medium equation is the presence of the
%logarithm term, which suggests that $v$ not only be bounded in
%magnitude, but also bounded away from $0$.

In Section~\ref{sec:unique} a result of Leon Simon~\cite{LS} is
used to prove Theorem~\ref{thm:unique}.

\section{Aronson-B\'{e}nilan Inequalities}\label{sec:bounds}

Let $\tilde u,v$ be solutions of Equations~\eqref{eq:secundus}
and~\eqref{eq:veq}, respectively.  The fundamental pointwise
inequalities upon which the other crucial bounds depend are
\begin{align}
t\frac{\D\tilde u}{\D t}&>-\frac{1}{2}\tilde u,\label{ineq:fund1}\\
(1-e^{-\tau})\frac{\D v}{\D\tau}&>-\frac{1}{2}v,\label{ineq:fund2}
\end{align}
and the integrated versions
\begin{align}
\tilde u(p,t_2)&>\sqrt{\frac{t_1}{t_2}}\tilde u(p,t_1),\label{ineq:intfund1}\\
v(p,\tau_2)&>e^{-\frac{\tau_2-\tau_1}{2}}
\sqrt{\frac{1-e^{-\tau_1}}{1-e^{-\tau_2}}}v(p,\tau_1),\label{ineq:intfund2}
\end{align}
for $t_2>t_1$ and $\tau_2>\tau_1$. To get these, we need only
assume that $\gamma$ is fixed and $\D f/\D\tau\geq 0$.

These are proved by an Aronson-B\'{e}nilan type argument similar to
that used for the porous medium equation~\cite{AB}.  To implement
this here, we define $w=1/\tilde u$ so that
Equation~\eqref{eq:secundus} becomes
\begin{equation}\label{eq:w}
    \frac{\D w}{\D
    t}=-\left(\Delta +f\right)w^{-1}.
\end{equation}
Defining now
\[
z\equiv t\frac{\D w}{\D t}-\frac{1}{2}w =-t\left(\Delta
+f\right)w^{-1}-\frac{1}{2}w,
\]
one finds
\begin{align}
  z'=-\left(\Delta +f\right)w^{-1}+t\left(\Delta
  +f\right)w^{-2}w'-f'tw^{-1}
  +\frac{1}{2}\left(\Delta +f\right)w^{-1},
\end{align}
where time differentiation has been denoted by a prime.  It is now
easily seen, using  Condition~\eqref{cond:2}, that $z$ satisfies
the linear parabolic differential inequality
\begin{equation}
z'\leq\left(\Delta+f\right)w^{-2}z.
\end{equation}
By the parabolic maximum principle, since $z$ is negative
initially, it must remain so.  Whence
\[
    t\frac{\D w}{\D t}<\frac{1}{2}w,
\]
and this inequality is equivalent to
Inequalities~\eqref{ineq:fund1}~and~\eqref{ineq:fund2}.

\section{Proof of Theorem~\ref{thm:bound}}\label{sec:global}
As indicated in the introduction, in this section
Theorem~\ref{thm:bound} is proved by using the maximum principle to
show that $\inf v$ must remain bounded. This establishes the result since,
also in this section, we obtain a Harnack inequality that bounds $\sup v$
in terms of $\inf v$.  The
latter is contained in
Proposition~\ref{prop:harnack}, whose proof is a direct
consequence of the weak Harnack inequalities of the next three
lemata, which are generalizations of results of Caffarelli and
Friedman~\cite{CF}.  The first of these establishes a lower bound
on $\inf v$ in terms of $\int_{\Sigma} v$, and the second and
third, in turn, use $\int_{\Sigma}v$ to bound $\sup v$ from above.

\begin{lemma}\label{lemma:one}
Assume $f,\, \D_{\tau}f\geq 0$,
and let $v$ be a solution of Equation~\eqref{eq:veq} with $v\geq\mu>0$.  Let $h>0$. There exist positive constants
$C,\tau_0$ with $C$ depending on $h$ and $\tau_0$ not depending on
$h$ such that
\[
   \int_{\Sigma} v(q,\tau_1)dV_q\leq C\left(\frac{1}{\mu}+\inf_{p\in\Sigma}v(p,\tau_1+h)\right),
\]
for any $q\in \Sigma$ and $\tau_1>\tau_0$.
\end{lemma}
\begin{remark}
The symbol $dV$ refers to the volume element with respect to $\gamma$, and $\Vol(\Sigma)$ will 
refer to the total volume of $\Sigma$ with respect to 
$\gamma$.  When it is clear from the context, $\int$ should be taken to mean $\int_{\Sigma}$, 
and if the volume element is omitted this should 
be taken as $dV$.  These remarks will continue to apply for the remainder of the text.  
\end{remark}

\begin{proof}
Let $G(p,q)$ be the positive Green's function such that
\[
     v(p,\tau)=-\int\Delta v(q,\tau)G(p,q)dV_q+\frac{1}{\Vol_{\gamma}(\Sigma)}\int
     v(q,\tau)dV_q;
\]
see~\cite{AU}. Using~\eqref{eq:veq} and rearranging, one has
\begin{align*}
    \frac{1}{\Vol(\Sigma)}\int v(q,\tau)dV_q
    &=\int\left(\frac{1}{v^2}\frac{\D v}{\D\tau}+\frac{1}{2v}-fv\right)G(p,q)dV_q+v(p,\tau)\\
    &\leq-\int\left(\frac{\D
    w}{\D\tau}-\frac{1}{2}w\right)G(p,q)dV_q+v(p,\tau),
\end{align*}
where  $w=1/v$. Multiplying by the integrating factor
$e^{-\frac{\tau}{2}}$, this becomes
\begin{equation}\label{ineq:green}
    e^{-\frac{\tau}{2}}\frac{1}{\Vol(\Sigma)}\int v(q,\tau)dV_q
    \leq-\int\frac{\D}{\D\tau}\left(e^{-\frac{\tau}{2}}w\right)G(p,q)dV_q+e^{-\frac{\tau}{2}}v(p,\tau).
\end{equation}
Assuming $\tau_1>\tau_0$, $\tau\in(\tau_1,\tau_1+h)$, and choosing
$\tau_0$ large enough that $t=1-e^{-\tau}>1/4$  for $\tau>\tau_0$, from~\eqref{ineq:intfund2} we have
\[
       \int v(q,\tau)dV_q\geq \frac{1}{2}e^{-\frac{(\tau-\tau_1)}{2}}\int v(q,\tau_1)dV_q
       \geq\frac{1}{2}e^{-\frac{(h)}{2}}\int v(q,\tau_1)dV_q
\]
and similarly
\[
    v(p,\tau)\leq 2e^{h/2}v(p,\tau_1+h).
\]
Using these two inequalities in~\eqref{ineq:green} we get
\begin{equation}
e^{-\frac{\tau}{2}}\frac{1}{\Vol(\Sigma)}\int v(q,\tau_1)dV_q
    \leq-2e^{h/2}\int\frac{\D}{\D\tau}\left(e^{-\frac{\tau}{2}}w\right)G(p,q)dV_q+
    e^{-\frac{\tau}{2}}4e^{h}v(p,\tau_1+h).
\end{equation}
Integrating over $(\tau_1,\tau_1+h)$ yields
$$\begin{array}{l}
\displaystyle\frac{e^{\frac{\tau_1}{2}}-e^{\frac{\tau_1+h}{2}}}{2}\frac{1}{\Vol(\Sigma)}\int
v(q,\tau_1)dV_q\\
    \displaystyle\leq  2e^{h/2}\int\left(w(q,\tau_1)e^{\frac{\tau_1}{2}}
    -w(q,\tau_1+h)e^{\frac{\tau_1+h}{2}}\right)G(p,q)dV_q\\
   \displaystyle+  \frac{e^{\frac{\tau_1}{2}}-e^{\frac{\tau_1+h}{2}}}{2}4e^{h}v(p,\tau_1+h).
\end{array}$$
Using now that $w\leq 1/\mu$, we obtain from this, finally
\[
     \frac{1}{\Vol(\Sigma)}\int
     v(q,\tau_1)dV_q\leq\frac{4e^{\frac{h}{2}}}{\mu\left(1-e^{\frac{h}{2}}\right)}\int
     G(p,q)dV_q+4e^{h}v(p,\tau_1+h).
\]
\end{proof}

\begin{lemma}\label{lemma:subsol}
Let $v$ be a solution of Equation~\eqref{eq:veq} on $[1,\infty)$
satisfying
\[
      v\geq\mu.
\]
Fix $\tau\in [1,\infty)$ and let
$M=\sup_{\Sigma}v(q,\tau),\,f^*=\sup_{\Sigma}f(q,\tau)$. Let $p\in\Sigma$, and 
let $r$ denote the geodesic distance from $p$. Then
given $\varepsilon>0$ there exists $r_0>0$ such that for $r<r_0$
there holds
\begin{equation}\label{ineq:smallharnack}
v(p,\tau)<
          \left(1+\varepsilon\right)\left(\frac{r^2}{2}\left(\frac{1}{2\mu}\frac{e^{-1}}{1-e^{-1}}+f^*M\right)
         +\frac{1}{|B_{r}(p)|}\int_{B_{r(p)}}v\right)
\end{equation}
for any $p\in\Sigma$.
\end{lemma}

\begin{proof}
For
\[
    H_0=\frac{1}{2}\left(\frac{1}{2\mu}\frac{e^{-1}}{1-e^{-1}}+f^*M\right)
\]
define $\phi=v+H_0r^2$.  Now, on a  neighborhood of $p$ let
$(r,\theta^1,\theta^2,...,\theta^{n-2})$ be geodesic polar
coordinates for $\gamma$ at $p$, for which we may write
$\gamma=dr^2+r^2h_{AB}d\theta^Ad\theta^B$.  Then (see~\cite{AU} p.
20) for a number $b$ such that $b^2$ bounds the sectional
curvature of $(\Sigma,\gamma)$ from above one has
   \[
      \Delta r^2=2(n-1)+2r\D_r\log |h|\geq 2(n-1)+2r\frac{\D}{\D
      r}\log\frac{\sin br}{r},
   \]
for $r_0$ small enough.  
And so, given $\delta>0$ we may choose $r_0$ small enough such
that $\Delta r^2\geq 2(n-1)-\delta$ whenever $r<r_0$. Using now
Equation~\eqref{eq:veq}, we have
\[
      \Delta\phi=\frac{1}{v^2}\frac{\D
      v}{\D\tau}+\frac{1}{2v}-fv+\left(2(n-1)-\delta\right)H_0.
\]
Hence by~\eqref{ineq:intfund2}
\begin{align*}
      \Delta\phi &>
      -\frac{1}{2v(1-e^{-\tau})}+\frac{1}{2v}-fv+\left(2(n-1)-\delta\right)H_0\\
      &=-\frac{1}{2v}\frac{e^{-\tau}}{1-e^{-\tau}}-fv+\left(2(n-1)-\delta\right)H_0\\
      &\geq
      -\frac{1}{2\mu}\frac{e^{-1}}{1-e^{-1}}-f^*M+\left(2(n-1)-\delta\right)H_0>0
\end{align*}
for $\delta<1$. Hence  $\phi$ is subharmonic, and  the lemma follows
from the mean value inequality for subharmonic functions  by
choosing $r_0$, perhaps, smaller still.
\end{proof}
We are now in a position to bound $\sup v$ in terms of an integral
of $v$.
\begin{lemma}\label{lemma:sup}
Let $v$ be a solution of Equation~\eqref{eq:veq} on $[1,\infty)$
satisfying $v\geq\mu$, and let $p,r$ be as in the previous lemma.  Given $\varepsilon>0$ there is an
$r_0$ such that for all $r\leq r_0$ one has
\[
     \sup_{p\in\Sigma}v(p,\tau)\leq
     2\left(1+\varepsilon\right)
    \left(\frac{r^2}{4\mu}\frac{e^{-1}}{1-e^{-1}}+\frac{1}{|B_{r}|}\int_{B_{r(p)}}v\right).
\]
\end{lemma}

\begin{proof}
For $\tau\in[1,\infty)$ fixed define $M=\sup_{\Sigma}v(p,\tau)$
and  let $p$ be such that $v(p,\tau)=M$.  Then by
Lemma~\ref{lemma:subsol}, assuming $r_0$ is small enough, one has
\[
M<
   \left(1+\varepsilon\right)\left(\frac{r^2}{2}\left(\frac{1}{2\mu}\frac{e^{-1}}{1-e^{-1}}
   +f^*M\right)+\frac{1}{|B_{r}|}\int_{B_{r(p)}}v\right).
\]
Thus as long as $r_0<\sqrt{1/(1+\varepsilon)f^*}$ there holds
\[
    M\leq
    2\left(1+\varepsilon\right)
    \left(\frac{r^2}{4\mu}\frac{e^{-1}}{1-e^{-1}}
    +\frac{1}{|B_{r}|}\int_{B_{r(p)}}v\right).
\]
\end{proof}
\begin{prop}\label{prop:harnack}
Assume $f,\, \D_{\tau}f\geq 0$,
and let $v$ be a solution of Equation~\eqref{eq:veq} with $v\geq\mu>0$. Let $h>0$. There exists  a constant $C$ independent
of $\tau$ but depending on $h$ such that
\[
   \sup_{p\in\Sigma}v(p,\tau)\leq
   C\left(\frac{1}{\mu}+\inf_{p\in\Sigma}v(p,\tau+h)\right).
\]
\end{prop}
\begin{proof}Apply Lemmata~\ref{lemma:one} and~\ref{lemma:sup}.\end{proof}
Theorem~\ref{thm:bound} can now be proved.
\begin{proof}[proof of Theorem~\ref{thm:bound}]
As previously remarked, the result follows directly from Proposition~\ref{prop:harnack} if
we can show that $\tilde v(\tau)=\inf_{\Sigma}v(p,\tau)$
remains bounded for all time.  In fact, there holds $\tilde
v<1/\sqrt{2\inf f}$. To see this, suppose instead that there is a
time $\tau_1$ at which $\tilde v(\tau_1)=1/\sqrt{2\inf f}$.  Then $v_*$
satisfying
\begin{align*}
    v_*'&=\inf f v_*^3+\inf f v_*^3-\frac{v_*}{2}\\
    v_*(\tau_1)&=\frac{1}{\sqrt{2\inf f}}
\end{align*}
is a subsolution of Equation~\eqref{eq:veq} for $\tau\geq\tau_1$ in the sense that
\[
   \frac{\D v_*}{\D\tau}\leq v_*^2\Delta v_*+fv_*^3-\frac{v_*}{2}
\]
for $\tau\geq\tau_1$.  
The parabolic maximum principle shows that $v\geq v_*$.  But $v_*$
blows up in finite time, and thus $v$ must also, which is a
contradiction to the definition of $v$.  Indeed, $v$ was taken to be a solution on $[\tau_0,\infty)$.  
\end{proof}
Finally, before leaving this section, we remark that any solution $v$ of
Equation~\eqref{eq:veq} satisfying the bounds of Theorem~\ref{thm:bound}
will in fact be uniformly bounded in $C^k$. That is, for every $\tau$ (large enough, of course) one
will have $\norm{v(\tau)}_{C^{k}(\Sigma)}\leq C$ for some constant only depending on $\mu,M$.
  The crucial
step towards doing this is to observe that on any finite
interval of the form $I=[0,T]$  such a solution $v$ will be uniformly bounded in
the parabolic analogue $H^{\alpha}_I$ of $C^{\alpha}$; for the
precise definition of $H^{k,\alpha}_I$ see~\cite{SW}. This
H\"{o}lder continuity follows from estimates originally due to Moser; see~\cite{lieberman}~Theorem
6.28. 
%\footnote{In order to obtain the H\"{o}lder regularity using
%Moser's techniques one must render the equation into divergence
%form.  In our case this gives gradient squared terms that threaten
%to violate the necessary structure conditions of, for
%instance,~\cite{lieberman}. However, careful checking shows that
%they do not. The reason for this is that in the case that we have
%good $L^{\infty}$ bounds on $u$, Moser's techniques can be
%implemented without putting the equation into divergence form,
%thereby avoiding the introduction of these dangerous gradient
%terms in the first place; see~\cite{SS}}. 
Afterwards, we may repeatedly apply
standard parabolic Schauder theory to get that $\norm{v}_{H^{k,\alpha}_{I}}\leq C$.
The desired bounds follow since, given
$\tau\in I$, one has $\norm{v(\tau)}_{C^{k,\alpha}(\Sigma)}\leq \norm{v}_{H_I^{k,\alpha}(\Sigma)}$.

\section{Proof of Theorem~\ref{thm:profile}}\label{sec:profile}
The  main ingredient used in the proof of
Theorem~\ref{thm:profile} is the functional
\[
    J(v)=\int_{\Sigma}|\nabla v|^2-fv^2+\log v,
\]
which is easily seen to be non-increasing in $\tau$~ by virtue of
Equation~\eqref{eq:veq} together with the condition $\D_{\tau}f\geq 0$:
\begin{equation}\label{ineq:lyap}
\frac{\D J}{\D\tau}=-2\int_{\Sigma}\frac{\D v}{\D\tau}\left(\Delta
v+fv-\frac{1}{2v}\right)-\int_{\Sigma}\frac{\D f
}{\D\tau}v^2\leq-2\int_{\Sigma}\frac{1}{v^2}\left|\frac{\D
v}{\D\tau}\right|^2
\end{equation}
The  hypothesis $\mu\leq v\leq M$  then establishes
\begin{align}
      J(v)&\geq-\sup f M^2+\log\mu\\
      \int_{\Sigma}|\nabla v|^2&\leq J(v_0)+\sup f M^2-\log\mu\label{ineq:intnablav}.
\end{align}
From the latter, we immediately obtain:
\begin{prop}\label{prop:L1conv}
Assuming the hypotheses of Theorem~\ref{thm:profile}, there exists
a sequence $\tau_i$ such that $v(\tau_i)\to \omega$ weakly in
$H^1$ and strongly in $L^1$.
\end{prop}
\begin{proof}
Following the argument of the introduction to this section, the
bound~\eqref{ineq:intnablav} shows that $v$ is bounded in $H^1$.
The result now follows from Rellich's theorem.  In the case $n=3$,
in which $\Sigma$ is 2-dimensional, one may apply Rellich's
theorem to $v$ as a function of $\Sigma \times \Sphere$, for
instance.
\end{proof}
We are now in a position to prove
Theorem~\ref{thm:profile}, for which it only remains to be shown that
$\omega$ is a stationary state and the convergence is actually in
$C^k$ for any $k$.
\begin{proof}[proof of theorem~\ref{thm:profile}]
Let  $v(\tau_i)$ be as in the conclusion of the preceding
proposition, fix $T>0$, and let $h<T$. Using the bound $v\leq M$
and~\eqref{ineq:lyap}, we get
\[
    \frac{1}{M^2}\int_{\Sigma}\left|\frac{\D
    v}{\D\tau}\right|^2\leq-\frac{1}{2}\frac{d}{d\tau}J(v).
\]
Integrating over $[\tau_i,\tau_i+h]$ and two applications of
H\"{o}lder's inequality yields
\[
    \norm{v(\tau_i+h)-v(\tau_i)}_{L^1}\leq
    C\sqrt{h}\sqrt{J(v(\tau_i))-J(v(\tau_i+h)))},
\]
for $C\geq\Vol_{\gamma}(\Sigma)$. Since we know that the
right hand side converges, we get that   $v(\tau_i+h)\to w$ in
$L^1$, uniformly for $h\in [0,T]$.

Now since $\mu\leq v\leq M$ it follows that $\mu\leq
\omega\leq M$ a.e.; hence
\[
     \int_{\Sigma} \left|\frac{1}{v(\tau_i+h)}-\frac{1}{\omega}\right|\leq
     \frac{1}{\mu^2}\int_{\Sigma} |\omega-v(\tau_i+h)|,
\]
and so $v^{-1}(\tau_i+h)\to {\omega}^{-1}$ in $L^1$ as well.

We are now in a position to prove that the limiting function
$\omega$ is a solution of the stationary equation.
Let $\psi$ be a $C^{\infty}$ function on $\Sigma$, and
$\varphi(\tau)$ a $C^{\infty}$ function compactly supported on
$[0,T]$, and put $v_i(p,\tau)=v(p,\tau_i+\tau)$. Then
\[
    \int_{0}^{T}\int_{\Sigma}\frac{\psi\varphi'}{v_i}
    =\int_{0}^{T}\int_{\Sigma}-\varphi\nabla\psi\cdot\nabla v_i+f\varphi\psi
    v_i-\frac{\varphi\psi}{v_i},
\]
and the convergence results from the previous paragraph show that
\[
     0=\int_{0}^{T}\int_{\Sigma}-\varphi\nabla\psi\cdot\nabla\omega+f\varphi\psi
    \omega-\frac{\varphi\psi}{\omega}.
\]
Hence $\omega$ is a weak solution of the stationary equation that
is essentially bounded above and below by positive constants. The
fact that it is a  $C^{\infty}$ solution follows from the
Sobolev embedding theorem and elliptic regularity.

To complete the proof, it only remains to show that the
convergence may also be taken in the sense of $C^k$. To do so,  we
look at the differences
 $\delta v_i=v_i-\omega$, which  verify
\[
      \delta v_i'=v _i^2\Delta\delta
      v_i+\left(fv_i^2+\frac{1}{2v_i\omega}\right)\delta
      v_i+(f_i-f)v_i,
\]
where we have put $f_i(\tau)=f(\tau+\tau_i)$.   Note that by the
remarks at the end of Section~\ref{sec:bounds}, the $v_i$ are
uniformly bounded in $H^{\alpha}_I$ for finite intervals
$I=[0,T]$. Thus, we may regard the previous equations as  linear
equations with uniformly H\"{o}lder continuous coefficients. By
the standard parabolic regularity theory we may conclude
     \[
      \norm{\nabla_{\gamma}^2\delta v_i}_{L^{q}(\Sigma\times\tilde I)}+\norm{ \D_{\tau}\delta v_i}_{L^{2}(\Sigma\times\tilde I)}\leq C
      \left(\norm{\delta v_i}_{L^q(\Sigma\times I)}+\norm{f_i-f}_{L^{q}(\Sigma\times I)}\right),
\]
for any $q>1$ and $\tilde I=[h_0,T], h_0>0$; see~\cite{lieberman}
p. 172. But since the  $v_i$  are uniformly bounded, convergence
of $\delta v_i$ to $0$ in $L^1$ implies convergence in $L^q$,
whereupon the previous inequality yields  convergence with respect to the parabolic analogue of $W^{2,q}$,
and the Sobolev embedding theorem implies $H^{\alpha}$ convergence
for some $0<\alpha<1$. This implies $H^{k,\alpha}$ convergence
since by the parabolic Schauder theory we have a bound of the form
\[
      \norm{\delta v_i}_{H^{k,\alpha}_{\tilde I}}\leq C
      \left(\norm{\delta v_i}_{C^0(\Sigma\times I)}+\norm{f_i-f}_{H^{\alpha}_I}\right).
\]
The theorem follows.
\end{proof}

\section{Proof of Theorem~\ref{thm:unique}}\label{sec:unique}
In this section it is proved that the stationary state $\omega$
found in the last section as the limit of certain sequences
$v(\tau_i)$ is not only unique, but in fact
$\lim_{\tau\to\infty}v(\tau)=\omega$, where the limit may be taken
in the sense of $C^{k}$ for any $k\in\mathbb N$. This will be proved
by using a result of Leon Simon~\cite{LS}. 

To describe Simon's result, consider the general case of a
parabolic equation
\begin{equation}\label{eq:upeq}
\frac{\D \nu}{\D\tau}=\EuScript M (\nu) + F, 
\end{equation}
where $\EuScript M$ is a second order elliptic differential operator
\[
\EuScript M:C^{\infty}(\Sigma)\to C^{\infty}(\Sigma),
\]
and $F$ is a smooth function on $\Sigma\times \R^+$, which is assumed to satisfy  
an exponential decay to be made more precise later.   

We assume furthermore that $\EuScript M$ is the gradient of an energy functional: there
exists
$\EuScript E:C^{\infty}(\Sigma)\to\R$ such that
\[
  <\EuScript M(\nu),\xi>_{L^2\left((\Sigma,\gamma)\right)}=-\frac{d}{ds}\EuScript
  E(\nu+s\xi)|_{s=0}
\]
for any $\xi,\nu\in C^{\infty}(\Sigma)$. The functional $\EuScript E$, 
in turn, is assumed to arise as the integral of an energy function: 
\[
   \EuScript E(\nu)=\int E(q,\nu,\nabla\nu)
\]
for some smooth $E:\times M\times \R\times T_pM\to \R$.  We
assume $E$ to be analytic, uniformly in $q$, as a function on
$\R\times T_qM$ in the sense that there exists $\beta$ such that
\begin{equation}\label{cond:analytic}
      E(q,z_0+s_1z,\overrightarrow p_0+s_2\overrightarrow p)=\Sigma_{|\alpha|\geq
      0}E_{\alpha}(q,z_0,z,\overrightarrow p_0,\overrightarrow p)s^{\alpha},\,\, s=(s_1,s_2),
\end{equation}
whenever $|z_0|,|z|,|\overrightarrow p_0|, |\overrightarrow p|<\beta$ and $|s|<1$;  
in addition, for these $z_0,z,p_0,p$ and for $j\geq 1$ we assume 
\begin{equation}
 \sup_{|s|<1}\left|\Sigma_{|\alpha|=j}E_{\alpha}(q,z_0,z,\overrightarrow p_0,\overrightarrow p)s^{\alpha}\right|\leq 1.
\end{equation}
We assume $E$ to be uniformly convex 
 in the sense that 
\begin{equation}\label{cond:convex}
 \frac{d^2}{ds^2}E(q,0,s\overrightarrow p)\big|_{s=0} \geq c|\overrightarrow p|^2
\end{equation}
for $c>0$ independent of $\tau, q$ and $\overrightarrow p$. 

The statement of the next theorem is contained in Theorem~2 of~\cite{LS}.  
\begin{theorem}\label{thm:Simonstheorem}
Let $\EuScript M$ and $F$ be as above, where in addition we assume $\EuScript M (0)=0$. Let $\nu$ be a $C^{\infty}$ 
solution of Equation~\eqref{eq:upeq} on $[0,\infty)$.  Let $l$ be
sufficiently large that $C^2(\Sigma)\subset W^{l-1,2}(\Sigma)$. Then there exists a $\delta>0$ such that if for some $T\geq 0$ one has
$\norm{v(\cdot,T)}_{W^{l-1,2}}<\delta$ and 
\begin{equation}\label{cond:F}
    \norm{F}_{W^{l-1,2}}+\norm{\D_{\tau}F}_{W^{l-1,2}}+\norm{\D^2_{\tau}F}_{L^2}\leq \delta e ^{T-\tau},\,\, \tau\geq T,
\end{equation}
then there holds  
\[
    \lim_{\tau\to\infty}\left(|\D_{\tau}\nu|_{C^2(\Sigma)}+|\nu-\nu_{\infty}|_{C^2(\Sigma)}\right)=0,
\]
where $\nu_{\infty}$ is a $C^{\infty}$ solution of the stationary equation $\EuScript M(\nu_{\infty})=0$.          
\end{theorem}

Thus, if there exists a sequence $\tau_m$ such that $\lim_{m\to\infty}\nu(\tau_m)\to 0$ in the sense of $C^2$, 
then it must be the case that $v_{\infty}=0$ and $\lim_{\tau\to\infty}\nu(\tau)=0$.  This will be the case for
us upon defining $\nu=\log(v/\omega)$.  But first, we must verify that $\nu$ satisfies an equation of the 
form~\eqref{eq:upeq}, for which all of the conditions listed above hold.    

The equation for $\nu$ is 
\begin{equation} \label{eq:neq}
      \frac{\D\nu}{\D\tau}=e^{\tilde\omega+\nu}\Delta e^{\tilde\omega+\nu} + f e^{\tilde\omega+\nu} -\frac{1}{2}, 
\end{equation}
where $\tilde\omega\equiv\log\omega$.  By defining 
\begin{equation}
    \EuScript M = e^{\tilde\omega+\nu}\Delta e^{\tilde\omega+\nu} + f_{t_1} e^{2\left(\tilde\omega+\nu\right)} -\frac{1}{2}
\end{equation}
and 
\begin{equation}
     F = \left(f-f_{t_1}\right)e^{2\left(\tilde\omega+\nu\right)}=  \left(f-f_{t_1}\right) v^2
\end{equation}
our equation takes the form~\eqref{eq:upeq} and $\EuScript M(0)=0$.  By the remarks at the end of Section~\ref{sec:bounds} we have that $v^2$ is 
bounded in $C^k(\Sigma)$ for every $k\in\mathbb N$, uniformly in $\tau$, and thus by the decay assumption on $f$, Condition~\eqref{cond:2'}, we see that 
$F$ satisfies~\eqref{cond:F}.  To continue, we note that our $\EuScript M$ is the gradient of the energy functional 
\[
     \EuScript E = \int E(q,\nu,\nabla\nu),
\]
where 
\[
    E(q,z,\overrightarrow p) = \frac{1}{2}\left(e^{2(\tilde\omega(q)+z)}
                               \left(\left|\overrightarrow p+\left(\nabla\tilde\omega\right)(q)\right|^2-f\right)+z\right).
\]
This can easily be checked to satisfy the analyticity and convexity assumptions~\eqref{cond:analytic}-\eqref{cond:convex}.  
Indeed, 
\[
       \frac{d^2}{ds^2}E(q,0,s\overrightarrow p) = e^{2\tilde\omega}|\overrightarrow p|^2\geq\mu^2 |\overrightarrow p|^2,
\]
and the analyticity is assumption is satisfied since the function $E$ is obtained as sums and products of linear, quadratic, 
and exponential functions in $z$ and $\overrightarrow p$. 

Now, by Theorem~\ref{thm:profile}, we have a sequence $\tau_i$ such that $\lim_{i\to\infty}\nu(\tau_i)=0$ in the sense of $C^2$.  
The remarks after 
Theorem~\ref{thm:Simonstheorem} imply that in fact $\lim_{\tau\to\infty}\nu(\tau)=0$.  Thus 
\[
    \lim_{\tau\to\infty} v = \lim_{\tau\to\infty} e^{\nu+\tilde\omega} = e^{\tilde\omega} = \omega.
\]
Using standard regularity theory one obtains convergence in $C^k(\Sigma)$ for all $k\in\mathbb N$.

\section{Appendix A: Blow-up solutions generated by the curve shortening flow}   \label{Appendix1}
Recall that a closed, simply connected, parameterized curve in the
plane $\gamma(\tau)$ is said to flow by curve shortening flow if
\[
    \frac{\D\gamma}{\D\tau}=k\vec n,
\]
where $k$ is the curvature and $\vec n$ is the inward pointing
unit normal.   In the analysis it is often more convenient to use
the support function and normal angle;  given a point $p$ on
$\gamma$ the normal angle is the angle between the position vector
and the normal vector, and the support function is defined by
$S(\theta)=\vec n\cdot \gamma$.  Note that in the case that
$\gamma$ is convex the normal angle gives a paramterization of
$\gamma$ on $[0,2\pi]$.  The support function $S$ and curvature $k$  
satisfy
\[
     S_{t}=-\frac{1}{S_{\theta\theta}+S}
\]
and
\[
      k_{t}=k^2(k_{\theta\theta}+k).
\]
Thus, in the case that the curve is strictly convex ($k>0$),
solutions of the curve shortening flow yield solutions of the
parabolic scalar curvature equation on
$\Sigma=\Sphere\times\Sphere$ with the product metric by taking
$r^2R-\kappa\equiv 1$ and $u(\theta_1,\theta_2)=\kappa(\theta_1)$.

Now, it is well known~\cite{GH} that under the curve shortening
flow,  the curve will shrink to a point in finite time, and more
specifically, the enclosed area $A(t)$ behaves according to
$A(t)=A(0)-2\pi t$.  To
study this shrinking more precisely one considers the normalized
curve
\[
      \gamma(\cdot,t)
      =\left(\frac{\pi}{A(t)}\right)^{\frac{1}{2}}\left(\gamma(\cdot,t)-\gamma(\cdot,\omega)\right),
\]
where $\omega=A(0)/(2\pi)$.  The curvature of this normalized
curve  is given by
\[
   \tilde k=\left(\frac{A(t)}{\pi}\right)^{\frac{1}{2}}k.
\]
Since it is well known from the work of Gage and Hamilton  that
the normalized curve converges to the unit circle as
$t\to\omega$~\cite{GH}, it follows that the normalized curvature
converges to $1$.  Hence, $u$ has the behavior claimed:
\[
     u=\frac{v}{\sqrt{A(0)-2\pi t}},
\]
where $v(\theta_1,\theta_2)=\tilde k(\theta_1)$ is uniformly
bounded  and in general can be assumed to vary in $\theta_1$
simply by assuming that the starting curve is different from the
circle.

\section{Appendix B: Derivation of the parabolic scalar curvature equation}  \label{Appendix2}

To derive Equation~\eqref{eq:primus}, let $N$ be the unit normal
vectorfield to the foliation leaves $\Sigma_r=\{r\}\times\Sigma$
so that the metric is written
\[
   g=u^2dr^2+r^2\gamma=N\otimes N+h,
\]
where  $h=r^2\gamma$.  For the calculation below it is important
to note that
\[
  N=\frac{\nabla r}{|\nabla r|}=u^{-1}\frac{\D}{\D r},
\]
and $u=|\nabla r|^{-1}$.  In addition, with $\nabla$ the covariant
derivative  compatible with $g$, it is convenient to define
$\nablash$ to be the covariant derivative induced by $\nabla$ on
$\Sigma_r$, which is compatible with $h$; and we shall also use
the notation
\[
   \Deltash=\Delta_h=r^{-2}\Delta_{\gamma}.
\]
Finally, let $\Pi$ denote the orthogonal projection of any
tensorfield onto $T\Sigma_r$ and recall that $\chi=\Pi\nabla N$
and $H=g^{ij}\chi_{ij}=h^{ij}\chi_{ij}=\nabla_iN^i$. We may now
begin the calculation.

Equation~\eqref{eq:primus} results from
\begin{equation}\label{eq:main0}
   \frac{\D H}{\D N}=-\frac{\Delta_h
   u}{u}-\frac 12 \left(H^2+|\chi|^2\right)+\frac12 \left(R_h-R\right),
\end{equation}
upon making the substitutions
$h=r^2\gamma,H=(ru)^{-1}\Hb,\chi=r\chib/u, N=u^{-1}\D_r$.
Equation~\eqref{eq:main0} in turn results from inserting the Gauss
equation
\[
    R_{NN}=\frac12 \left(R-R_h\right)+\frac12 \left(H^2-|\chi|^2\right)
\]
into
\begin{equation}\label{eq:2ndvar}
\frac{\D H}{\D N}=-\frac{\Delta_h u}{u}-|\chi|^2 -R_{NN},
\end{equation}
where the latter is a key formula used in the general second variation of area. For the sake of completeness, we 
derive~\eqref{eq:2ndvar} and verify the relation of $\Hb,\chib$ to
$H,\chi$, which was stated in the introduction. The latter  shall be done first.

From Killing's formula and the definition of $\chi$ one has
\[
   \chi=\frac{\Pi}{2}\Lie_Ng=\frac{\Pi}{2}\left(\Lie_NN\otimes N+N\otimes\Lie
   N+\Lie_Nh\right)
   =\frac{1}{2}\Lie_Nh.
\]
Hence
\begin{align*}
    \chi=\frac12\Lie_Nh&=\frac12\Lie_{u^{-1}\D_r}h=\frac12\left(u^{-1}\Lie_{\D_r}h+\nabla
    u^{-1}\otimes h(\D_r,\cdot)+h(\cdot,\D_r)\otimes\nabla
    u^{-1}\right)\\
    &=\frac12u^{-1}\Lie_{\D_r}h=u^{-1}\frac12\frac{\D h_{AB}}{\D
    r}d\theta^Ad\theta^B,
\end{align*}
where in the last step we have used the coordinates
$(\theta^1,\theta^2,...,\theta^{n-1})$ from the introduction.  But
since $h=r^2\gamma$, we get the expression for $\chi$
\[
      \chi=ru^{-1}\left(\gamma_{AB}+\frac{1}{2}r\frac{\D\gamma_{AB}}{\D
      r}\right),
\]
and thus $\chi=ru^{-1}\chib$.  The relation of $H$ to $\Hb$ then
follows by contraction with $h^{AB}=r^{-2}\gamma^{AB}$.

We now derive Equation~\eqref{eq:2ndvar}:
\begin{align*}
\frac{\D H}{\D
N}&=N^k\nabla_k\nabla_iN^i=N^k\nabla_i\nabla_kN^i-R_{kil}^{\,\,\,\,\,\,\,i}N^kN^l\\
  &=\nabla_i\left(N^k\nabla_kN^i\right)-\nabla_iN^k\nabla_kN^i-R_{NN}\\
  &=\nabla\cdot\left(\nabla_NN\right)-|\chi|^2-R_{NN}
\end{align*}
The term $\nabla_N N$ may be calculated in terms of $u$ and $\nablash u$ as follows:
\begin{align*}
 \nabla_NN_j&=N^k\nabla_k\left(\frac{\nabla_j\,r}{|\nabla
 r|}\right)=N^k\nabla_k\left(u\nabla_j\,r\right)=\nabla_Nu\nabla_j\,
 r+uN^k\nabla_k\nabla_j\,r\\
 &=\frac{\nabla_N u}{u}N_j+uN^k\nabla_j\nabla_k\,r=\frac{\nabla_N
 u}{u}N_j+u\nabla_j\left(N^k\nabla_k\,
 r\right)-u\nabla_jN^k\nabla_k\,r\\
 &=\frac{\nabla_N u}{u}N_j+u\nabla_j\left(u^{-1}\right)=\frac{\nabla_N
 u}{u}N_j-\frac{\nabla_j\,u}{u}=-\frac{\nablash_j\,u}{u}
\end{align*}

Finally, the derivation of~\eqref{eq:2ndvar} is completed upon showing that
$\nabla\cdot\left(\nablash u/u\right)=\Deltash u/u$:
\begin{align*}
 \nabla\cdot\left(\nablash u/u\right)&=g^{ij}\nabla_i\left(\frac{\nablash_ju}{u}\right)
                                      =\left(h^{ij}+N^iN^j\right)\nabla_i\left(\frac{\nablash_ju}{u}\right)\\
                                     &=h^{ij}\nablash_i\left(\frac{\nablash_ju}{u}\right)+N^j\nabla_N\left(\frac{\nablash_ju}{u}\right)\\
                                     &=\left(\frac{\Deltash u}{u}-\frac{|\nablash u|^2}{u^2}\right)
                                       +\nabla_N\left(N^j\frac{\nablash_j u}{u}\right)-\nabla_NN^j\left(\frac{\nablash_ju}{u}\right)\\
                                     &=\left(\frac{\Deltash u}{u}-\frac{|\nablash u|^2}{u^2}\right)
                                       +0+\frac{|\nablash u|^2}{u^2}=\frac{\Deltash u}{u}.
\end{align*}

\bibliographystyle{amsplain}

\end{document}